# A reference-searching–based algorithm for large-scale data envelopment analysis computation

Wen-Chih Chen*

Department of Industrial Engineering and Management

National Chiao Tung University, Taiwan

## Abstract

Data envelopment analysis (DEA) is a linear program (LP)-based method used to determine the efficiency of a decision making unit (DMU), which transforms inputs to outputs, by peer comparison. This paper presents a new computation algorithm to determine DEA efficiencies by solving small-size LPs instead of a full-size LP. The concept is based on searching the corresponding references, which is a subset of the efficient DMUs with numbers no greater than the dimension (number of inputs and outputs). The results of empirical case studies show that the proposed algorithm computes 3 times faster than the current state of the art for large-scale, high-dimension, and high-density (percentage of efficient DMUs) cases. It is flexible enough to compute the efficiencies of a subset of the full data set without setup costs, and it can also serve as a sub-procedure for other algorithms.

* Corresponding author. 1001 Ta Hsueh Road, Hsinchu, Taiwan. Email: wenchih@faculty.nctu.edu.tw

## 1. Introduction

Data envelopment analysis (DEA), originated by Farrell (1957) and popularized by Charnes et al. (1978), measures the efficiency of a firm, or in general a decision making unit (DMU) by peer comparison. DEA considers all DMUs in the peer transforming inputs to outputs under similar environments, identifies the ideal transformation as the "efficient frontier", and measures the "gap" of a DMU's transformation outcome to the frontier as the efficiency. From the computational geometry perspective, DEA considers a set of finite data points, each corresponding to a DMU's input-output transformation result, in a high dimensional space. The tightest polyhedral containing all the points is identified and its relative location with respect to a particular data point is measured based on a pre-specified definition. A DMU is termed efficient if its input-output data is on the efficient frontier; otherwise, it is inefficient.

The implementation of determining a DMU's efficiency solves a linear program (LP) as the analytical and computational technique. A standard procedure is to formulate and solve an LP for each DMU to obtain all efficiencies of the peers. Each LP has a dense matrix generated by the full data set (i.e. each DMU contributes its input-output data as a column), and the number of data points and the number of inputs and outputs (termed dimension) determine the size of the matrix. An increase in the size of the data set requires a higher number of larger full-size LPs to solve. Although an LP is solvable in polynomial time, solving a higher number of full-size LPs becomes computationally burdensome when the data set is massive (cf. Barr and Durchholz, 1997 and Dula, 2011).

Much of the published literature on DEA computation for large-scale data sets focuses on reducing the size of individual LPs. The idea is generally based on the observation that only efficient DMUs affect the computation results, i.e. inefficient DMUs will not affect the computation results. For example, Dulá and Thrall (2001) and Dulá (2011) suggested identifying all of the extreme efficient[1] DMUs and using them to determine the efficiency of the other DMUs. Barr and Durchholz (1997) divided the DMUs into several subsets, solved the efficiencies of all of the DMUs within each subset, and removed the corresponding inefficient DMUs to reduce the size of individual LP solving. Sueyoshi and Chang (1989), Ali (1993), and Chen and Ali (2002) proposed preprocessors to reduce computation without the need to solve the LPs (see Dulá and López (2009) for a comprehensive review).

This paper develops a new algorithm for computing DEA efficiency from the computational geometry perspective of the input-output space. The algorithm searches the "references," which is a small critical subset of efficient DMUs to measure the efficiency with a size no larger than the dimension. The references compose a piece of the frontier (a facet) corresponding to the data point of a DMU, which then allows measuring the efficiency as the relative location between them. The searching process proceeds to solve several small LPs until the true references are included in the LP model.

The idea and framework of reference searching was discussed in Chen and Lai (2017). Focusing on that it is possible to compute a large-scale DEA problem by solving size-limited LPs, they proposed to switch several columns (corresponding to DMUs) in and out of the fixed-size LPs until the LP included the critical DMUs necessary to compute the efficiency. The discussion in Chen and Lai, however, is general and broad in computational performance; they did not discuss how to accelerate DEA computation or address the detailed mechanisms of the framework and the technical challenges encountered, such as multiple dual optima. Whereas the proposed algorithm in this paper presents very detailed theoretical analyses, refines and integrates the algorithms developed in Chen and Lai (2017) and represents the new direction in DEA computation.

[1] The term, extreme efficient, is defined in Charnes et al. (1991).

The proposed algorithm significantly enhances large-scale DEA computation, and it is particularly insensitive to dimension and density (percentage of efficient DMUs). The results of the empirical studies in Section 5 below show that the proposed algorithm can compute 2.5 times faster than the state of the art in large-scale, high-dimension, and high-density cases. Notably, it has been designed to determine the efficiency of one DMU, i.e. the algorithm needs to be processed once for each DMU to determine the efficiencies of all of them. This unique feature provides the flexibility to adopt our algorithm. For example, the computation load is linear to the number of DMUs without a fixed considerable "setup" cost. The proposed algorithm can also serve as a sub-procedure for other algorithms presented in the literature. Moreover, by solving a new LP formulation, the proposed algorithm can address the challenge of degeneracy, which leads to multiple multiplier weights in DEA (cf. Ali, 1994; Barr and Durchholz, 1997; Cooper et al., 2007).

The remainder of this paper is organized as follows. Section 2 introduces the notation, assumptions, and underlying concepts. Section 3 develops the theoretical foundation, and Section 4 presents the proposed algorithm. Section 5 discusses empirical case studies and benchmarks the proposed algorithm's performance with state of the art. Section 6 concludes and offers suggestions for future research.

## 2. Notation and assumptions

In this paper we consider a DEA study with a set of DMUs, $\mathcal{D}$, and the associated input-output data recording the outcomes of how the DMUs transform $m$ inputs to $n$ outputs. We denote DMU $r$'s input-output data as $(x_r, y_r) \in \Re^{m+n}$, where $x_r$ and $y_r$ are the input and output value vectors, respectively. We assume that the input-output values are strictly positive, i.e. $(x_r, y_r) \in \Re^{m+n}_{++}$ and all input-output vectors are "reduced", i.e. distinct and not proportionate to each other.

Modeling the input-output transformation and defining an efficiency measure are key components of DEA. DEA characterizes the feasible input-output transformation by the *production possibility set* and approximates it by the data associated with $\mathcal{D}$, based on economic assumptions such as convexity and free disposability. A popular approximation is:

$$\mathcal{T} = \left\{(x, y): x \geq \sum_{r\in\mathcal{D}} \lambda_r x_r\,; y \leq \sum_{r\in\mathcal{D}} \lambda_r y_r\,; \sum_{r\in\mathcal{D}} \lambda_r = 1; \lambda_r \geq 0, \forall r \in \mathcal{D}\right\}.$$

$\mathcal{T}$ is a $m+n$ dimensional polyhedral set defined on the data. It is a convex hull containing the data points and recedes along the direction with higher input values and lower output values. Farrell (1957) proposed to measure the efficiency as the maximum proportionate extraction of inputs or expansion of outputs. It determines how much of the proportionate extraction of inputs or expansion of outputs was required to move the data point of a DMU to the boundary of the production possibility set. Farrell's measure requires no *a priori* information regarding prices or weights of inputs and outputs, and is termed the *radial efficiency*.

Without loss of generality, we solve the following LP to determine the radial efficiency of $k \in \mathcal{D}$ (Banker et al., 1984),

$$\begin{aligned} e_k^* = \min \quad & \theta \qquad (1) \\ \text{s.t.} \quad & \sum_{r\in\mathcal{D}} \lambda_r x_r \leq \theta x_k; \\ & \sum_{r\in\mathcal{D}} \lambda_r y_r \geq y_k; \\ & \sum_{r\in\mathcal{D}} \lambda_r = 1; \\ & \lambda_r \geq 0, \forall r \in \mathcal{D}, \text{and } \theta \text{ is unrestricted.} \end{aligned}$$

The variables $\lambda_r$'s and $\theta$ have different roles in the computational geometry. $\lambda_r$'s define the linear combinations of the data points and relate to the production possibility set. The constraints of Problem (1) specify the feasible input-output transformation approximated by the data of DMU set $\mathcal{D}$. Any input-output vector $(\theta x_k, y_k)$ satisfying the constraints of problem (1) represents a feasible transformation and is in the production possibility set.

Variable $\theta$ and the objective function specify how $(x_k, y_k)$ is improved based on Farrell's input-oriented radial measure. Denoting $\theta^o$ and $\lambda_r^o$ as the optimal solutions of Problem (1), $e_k^* = \theta^o$ is the efficiency of DMU $k$. $(\theta^o x_k, y_k)$, which associates with the right side of the first two types of constraints, is the improvement limit for $k$, and is the projected point on the boundary of the production possibility set. Comparing to the data associated with set $\mathcal{D}$, $(x_k, y_k)$, we conclude that $(x_k, y_k)$ can be improved by reducing the inputs to $100 \times \theta^o\%$ of its current level while maintaining the current output level. DMU $k$ is efficient if $e_k^* = \theta^o = 1$ since there is no room to reduce the inputs proportionately; otherwise, DMU $k$ is inefficient.

$(\sum_{r\in\mathcal{D}} \lambda_r^o x_r, \sum_{r\in\mathcal{D}} \lambda_r^o y_r)$, which associates with the left side of the first two types of constraints, is the benchmark composited from the peer for $k$'s improvement. It is not necessary for $(\theta^o x_k, y_k)$ and $(\sum_{r\in\mathcal{D}} \lambda_r^o x_r, \sum_{r\in\mathcal{D}} \lambda_r^o y_r)$ to be identical due to the existing slacks in the constraints. Thus, an efficient DMU classified by $e_k^*$ may not be Pareto-Koopmans efficient, and is further classified as *weakly efficient*. In contrast, a Pareto-Koopmans efficient DMU is termed *strongly efficient*. The boundary of the production possibility set (or the frontier) can also be further classified by the types (weakly or strongly efficient) of the points on it. We focus on determining the DEA efficiency, not the classification, but our algorithm applies to different classifications of efficient DMUs.

$(\sum_{r\in\mathcal{D}} \lambda_r^o x_r, \sum_{r\in\mathcal{D}} \lambda_r^o y_r)$ is only influenced by points with $\lambda_r^o > 0$. We term $\{r \in \mathcal{D}: \lambda_r^o > 0\}$ the "reference points" or simply the "references" of $k$. Problem (1) has $|\mathcal{D}| + 1$ variables and $m + n + 1$ constraints. Typically, we have $|\mathcal{D}| \gg m + n$. The optimal solutions of Problem (1) have at most $m + n$ $\lambda$'s in the basis that are possibly non-zero. The last slot in the basis is always reserved for $\theta$. Therefore, regardless of problem size, at most $m + n$ DMUs in $\mathcal{D}$ are necessary for solving Problem (1), and they are references for DMU $k$.

The dual to Problem (1) is

$$\begin{aligned} e_k^* = \max \quad & v y_k + w \\ \text{s.t.} \quad & u x_k = 1; \\ & -u x_r + v y_r + w \le 0 \qquad \forall r \in \mathcal{D}; \\ & u, v \ge 0,\ w \text{ is unrestricted.} \end{aligned} \tag{2}$$

Here, $u \in \Re_+^m$, $v \in \Re_+^n$, and $w \in \Re$ are the dual variables corresponding to the constraints of Problem (1) from the top to the bottom. For notational simplicity, we denote the inner product of two vectors $u$ and $x_r$ as $ux_r$, and apply this notation to other vectors in this paper.

In $\Re^{m+n}$, any $(u,v,w)$ defines a hyperplane

$$\mathcal{H}^{(u,v,w)} = \{(x,y) \in \Re^{m+n} | -ux + vy + w = 0\}$$

and two half-spaces of interest: $\mathcal{H}_{++}^{(u,v,w)} = \{(x,y) \in \Re^{m+n} | -ux + vy + w > 0\}$ and $\mathcal{H}_{-}^{(u,v,w)} = \{(x,y) \in \Re^{m+n} | -ux + vy + w \leq 0\}$. While the sign of $-ux + vy + w$ positions $(x,y) \in \Re^{m+n}$ in two half-spaces, $\mathcal{H}_{-}^{(u,v,w)}$ and $\mathcal{H}_{++}^{(u,v,w)}$, the absolute value of $-ux + vy + w$ represents the "distance" of $(x,y)$ to the hyperplane $\mathcal{H}^{(u,v,w)}$. $-ux + vy + w < -ux' + vy' + w \leq 0$ implies that $(x',y') \in \Re^{m+n}$ is "closer" to $\mathcal{H}^{(u,v,w)}$ than $(x,y)$. Similarly, $-ux + vy + w > -ux' + vy' + w > 0$ implies that $(x,y)$ is "farther" from $\mathcal{H}^{(u,v,w)}$ than $(x',y')$.

Along this line, Problem (2) searches for $(u,v,w)$ such that $\mathcal{H}_{-}^{(u,v,w)}$ contains all points associated with $\mathcal{D}$ and $(x_k,y_k)$ is as close as possible to $\mathcal{H}^{(u,v,w)}$. $\left\{r \in \mathcal{D} \middle| (x_r,y_r) \in \mathcal{H}_{++}^{(u,v,w)}\right\} \neq \emptyset$ indicates that $(u,v,w)$ is infeasible to Problem (2). Shifting the hyperplane $\mathcal{H}^{(u,v,w)}$ to move the points from one side of the hyperplane ($\mathcal{H}_{++}^{(u,v,w)}$) to the other side ($\mathcal{H}_{-}^{(u,v,w)}$) provides a feasible solution to Problem (2). Clearly, when the point in $\mathcal{H}_{++}^{(u,v,w)}$ that is farthest from the hyperplane $\mathcal{H}^{(u,v,w)}$ switches its location with respect to the hyperplane, all points are in the same half-space, giving a feasible solution.

We denote the optimal solution of Problem (2) as $(u^o, v^o, w^o)$. At least one $r \in \mathcal{D}$ must exist such that $-u^o x_r + v^o y_r + w^o = 0$; otherwise, a better solution can be obtained. $\mathcal{H}_{-}^{(u^o,v^o,w^o)}$ is the "smallest" half-space containing all data points associated with $\mathcal{D}$; $\mathcal{H}^{(u^o,v^o,w^o)}$ is tight to the data set and termed the *supporting hyperplane,* or simply the *hyperplane,* corresponding to $k$. $-u^o x_k + v^o y_k + w^o = 0$ indicates that $(x_k,y_k)$ is on the hyperplane $\mathcal{H}^{(u^o,v^o,w^o)}$; it also follows $v^o y_k + w^o = 1$ and indicates that $e_k^* = 1$. $v^o y_k + w^o = e_k^* < 1$ implies $-u^o x_k + v^o y_k + w^o = -1 + e_k^* < 0$ (i.e. a point associated with an inefficient DMU is in the interior of $\mathcal{H}_{-}^{(u^o,v^o,w^o)}$). $(e_k^* x_k, y_k)$ gives $-u^o(e_k^* x_k) + v^o y_k + w^o = 0$ since $v^o y_k + w^o = e_k^*$ and $-u^o x_k = 1$, and is on the hyperplane. Therefore, $e_k^*$ provides information on the relative location between $(x_k, y_k)$ and $\mathcal{H}^{(u^o,v^o,w^o)}$.

**3. Theoretical development**

We base our theoretical development on an LP with respect to $\mathcal{D}$, $k \in \mathcal{D}$ and $M \geq 1$:

$$\theta_k^{\mathcal{D}} = \min \quad \theta \qquad \text{P}(\mathcal{D})$$

$$\text{s.t.} \quad \lambda_k M x_k + \sum_{r\in\mathcal{D}\backslash\{k\}} \lambda_r x_r \leq \theta M x_k;$$

$$\lambda_k y_k + \sum_{r\in\mathcal{D}\backslash\{k\}} \lambda_r y_r \geq y_k;$$

$$\lambda_k + \sum_{r\in\mathcal{D}\backslash\{k\}} \lambda_r = 1;$$

$$\lambda_k, \lambda_r \geq 0, \forall r \in \mathcal{D}\backslash\{k\}, \text{and } \theta \text{ is unrestricted.}$$

Problem $\text{P}(\mathcal{D})$ is identical to Problem (1) except that $(x_k, y_k)$ is replaced by $(Mx_k, y_k)$ (i.e. it determines the efficiency of $(Mx_k, y_k)$, not $(x_k, y_k)$). We note that Problem $\text{P}(\mathcal{D})$ is always feasible since $\theta = 1$, $\lambda_k = 1$, and $\lambda_r = 0, \forall r \in \mathcal{D}\backslash\{k\}$ is a feasible solution. The dual to Problem $\text{P}(\mathcal{D})$ is

$$\theta_k^{\mathcal{D}} = \max \quad v y_k + w \qquad \text{D}(\mathcal{D})$$

$$\text{s.t.} \quad u(Mx_k) = 1;$$

$$-u x_r + v y_r + w \leq 0 \qquad \forall r \in \mathcal{D}\backslash\{k\};$$

$$-u(Mx_k) + v y_k + w \leq 0;$$

$$u, v \geq 0, \ w \text{ is unrestricted.}$$

We denote $\theta^*$ and $\lambda_r^*$ as the optimal solutions of Problem $\text{P}(\mathcal{D})$ and $(u^*, v^*, w^*)$ as the optimal solution to Problem $\text{D}(\mathcal{D})$ hereafter. It follows that

$$u^*(Mx_k) = 1; \tag{3}$$

$$-u^* x_r + v^* y_r + w^* \leq 0, \forall r \in \mathcal{D}\backslash\{k\}; \tag{4}$$

$$-u^*(Mx_k) + v^* y_k + w^* \leq 0. \tag{5}$$

Eq. (3) also implies that $u^*$ cannot be a zero vector.

The first result shows how $e_k^*$ and $\theta_k^{\mathcal{D}}$ relate when $(x_k, y_k) \in H_-^{(u^*,v^*,w^*)}$.

**RESULT 1**: If $-u^*x_k + v^*y_k + w^* \leq 0$, then $e_k^* = M\theta_k^{\mathcal{D}}$.

**Proof:**

(i) Eq. (3) yields $(Mu^*)x_k = 1$. Eq. (4) and $-u^*x_k + v^*y_k + w^* \leq 0$ imply $-u^*x_r + v^*y_r + w^* \leq 0, \forall r \in \mathcal{D}$ and $-(Mu^*)x_r + (Mv^*)y_r + (Mw^*) \leq 0, \forall r \in \mathcal{D}$. Therefore, $(Mu^*, Mv^*, Mw^*)$ is feasible to Problem (2) and $Mv^*y_k + Mw^* = M(v^*y_k + w^*) \leq e_k^*$ (i.e. $M\theta_k^{\mathcal{D}} \leq e_k^*$).

(ii) Similarly, given $(u^o, v^o, w^o)$ as the optimal solution for Problem (2), we have $v^oy_k + w^o = e_k^*$, $u^ox_k = 1$, $-u^ox_r + v^oy_r + w^o \leq 0, \forall r \in \mathcal{D}$, and $u^o, v^o \geq 0$. $M \geq 1$ and $-u^ox_r + v^oy_r + w^o \leq 0$ follow $-\frac{u^o}{M}x_r + \frac{v^o}{M}y_r + \frac{w^o}{M} \leq 0$ and $-\frac{u^o}{M}(Mx_r) + \frac{v^o}{M}y_r + \frac{w^o}{M} \leq -\frac{u^o}{M}x_r + \frac{v^o}{M}y_r + \frac{w^o}{M} \leq 0$ $\forall r \in \mathcal{D}$. Together with $\frac{u^o}{M}(Mx_k) = 1$, since $u^ox_k = 1$, $\left(\frac{u^o}{M}, \frac{u^o}{M}, \frac{w^o}{M}\right)$ is feasible to Problem D($\mathcal{D}$) and $\frac{u^o}{M}y_k + \frac{w^o}{M} \leq \theta_k^{\mathcal{D}}$, i.e., $e_k^* \leq M\theta_k^{\mathcal{D}}$.

Combining (i) and (ii) gives $M\theta_k^{\mathcal{D}} = e_k^*$. □

P($\mathcal{D}$) is identical to Problem (1) when $M = 1$ and $\theta_k^{\mathcal{D}} = e_k^*$. However, $-u^*x_k + v^*y_k + w_0^* \leq 0$ is unnecessary when $M > 1$. For example, we have $-u^*x_k + v^*y_k + w^* > 0$ when $-u^*(Mx_k) + v^*y_k + w^* = 0$ ($\theta_k^{\mathcal{D}} = 1$) including both strongly and weakly efficient DMUs. The following result shows the relation between P($\mathcal{D}$) and (1) for $(x_k, y_k) \in H_{++}^{(u^*,v^*,w^*)}$.

**RESULT 2**: $-u^*x_k + v^*y_k + w^* > 0$ imply $e_k^* = 1$.

**Proof:**

When $M = 1$, we have $\theta_k^{\mathcal{D}} = e_k^*$ and $-u^*x_k + v^*y_k + w^* \leq 0$, and we need only discuss the cases with $M > 1$. Denoting $\tau = u^*y_k + w^*$, and together with (3) gives

$$-\tau Mu^*x_k + v^*y_k + w^* = 0. \tag{6}$$

$-u^*x_k + v^*y_k + w^* > 0$ and Eq. (3) gives $\tau M > 1$.

We have $-(\tau Mu^*)x_r + v^*y_r + w^* \leq 0, \forall r \in \mathcal{D}\backslash\{k\}$ because of (4) and $\tau M > 1$, and it follows $-\frac{(\tau Mu^*)x_r}{\tau} + \frac{v^*y_r}{\tau} + \frac{w^*}{\tau} = -(Mu^*)x_r + \frac{v^*}{\tau}y_r + \frac{w^*}{\tau} \leq 0, \forall r \in \mathcal{D}\backslash\{k\}$. Eq. (6) gives $(-\tau Mu^*x_k + v^*y_k + w^*)\frac{1}{\tau} = -(Mu^*)x_k + \frac{v^*}{\tau}y_k + \frac{w^*}{\tau} = 0$. Eq. (3) gives $(Mu^*)x_k = 1$.

Therefore, $\left(Mu^*, \frac{v^*}{\tau}, \frac{w^*}{\tau}\right)$ is feasible to Problem (2). $(Mu^*)x_k = 1$ also gives $\frac{u^*}{\tau}y_k + \frac{w^*}{\tau} = 1$ based on (6). Therefore, $e_k^* = 1$, and $\left(Mu^*, \frac{v^*}{\tau}, \frac{w^*}{\tau}\right)$ is optimal to Problem (2). □

Now we consider the problem where $\mathcal{S} \subseteq \mathcal{D}$ and $k \in \mathcal{S}$

$$\theta_k^{\mathcal{S}} = \min \quad \theta \qquad \text{P}(\mathcal{S})$$

$$\text{s.t.} \quad \lambda_k M x_k + \sum_{r \in \mathcal{S}\setminus\{k\}} \lambda_r x_r \le \theta M x_k;$$

$$\lambda_k y_k + \sum_{r \in \mathcal{S}\setminus\{k\}} \lambda_r y_r \ge y_k;$$

$$\lambda_k + \sum_{r \in \mathcal{S}\setminus\{k\}} \lambda_r = 1;$$

$$\lambda_k, \lambda_r \ge 0, \forall r \in \mathcal{S}\setminus\{k\}, \text{and } \theta \text{ is unrestricted.}$$

Apparently, $\theta_k^{\mathcal{S}} \ge \theta_k^{\mathcal{D}}$ since $\mathcal{S} \subseteq \mathcal{D}$, but an $\mathcal{S} \subseteq \mathcal{D}$, named $\mathcal{S}^*$, must exist such that $\theta_k^{\mathcal{S}^*} = \theta_k^{\mathcal{D}}$. Properly selecting the DMUs gives $\mathcal{S}^*$, which allows us to determine the efficiency with a massive data set by solving small-size LPs constructed by a sample with a small number of data points. We denote $\{\hat{\theta}; \hat{\lambda}_r, r \in \mathcal{S}\}$ as the optimal solution to Problem P($\mathcal{S}$) and $(\hat{u}, \hat{v}, \hat{w})$ as its optimal dual solution. The following result provides the sufficient condition for $\theta_k^{\mathcal{S}} = \theta_k^{\mathcal{D}}$.

**RESULT 3:** $\theta_k^{\mathcal{S}} = \theta_k^{\mathcal{D}}$ if the following inequalities are satisfied:

$$-\hat{u}x_t + \hat{v}y_t + \hat{w} \le 0, \forall t \in \mathcal{D}\setminus\mathcal{S}. \qquad (7)$$

**Proof:**

We propose $\{\hat{\theta}; \hat{\lambda}_r, r \in \mathcal{S}; \hat{\lambda}_t = 0, t \in \mathcal{D}\setminus\mathcal{S}\}$ and $(\hat{u}, \hat{v}, \hat{w})$ as the optimal solution to P($\mathcal{D}$) and its dual, respectively. We check the proposed solution with the Karush-Kuhn-Tucker (KKT) optimality condition corresponding to P($\mathcal{D}$) that satisfies the following constraints.

Primal feasibility:

$$\sum_{t \in \mathcal{D}\setminus\mathcal{S}} \hat{\lambda}_t x_t + \sum_{r \in \mathcal{S}\setminus\{k\}} \hat{\lambda}_r x_r + \hat{\lambda}_k (M x_k) \le \hat{\theta}(M x_k),$$

$$\sum_{t \in \mathcal{D}\setminus\mathcal{S}} \hat{\lambda}_t y_t + \sum_{r \in \mathcal{S}\setminus\{k\}} \hat{\lambda}_r y_r + \hat{\lambda}_k y_k \ge y_k,$$

$$\sum_{t\in\mathcal{D}\setminus\mathcal{S}} \hat{\lambda}_t + \sum_{r\in\mathcal{S}\setminus\{k\}} \hat{\lambda}_r + \hat{\lambda}_k = 1,$$

$\hat{\lambda}_t \geq 0 \ \ \forall t \in \mathcal{D}\setminus\mathcal{S}$ and $\hat{\lambda}_r \geq 0 \ \ \forall r \in \mathcal{S}$,

Dual feasibility:

$\hat{u}(Mx_k) = 1,$

$-\hat{u}x_r + \hat{v}y_r + \hat{w} \leq 0 \quad \forall r \in \mathcal{S}\setminus\{k\},$

$-\hat{u}(Mx_k) + \hat{v}y_k + \hat{w} \leq 0,$

$-\hat{u}x_t + \hat{v}y_t + \hat{w} \leq 0 \quad \forall t \in \mathcal{D}\setminus\mathcal{S},$

$\hat{u} \geq 0$ and $\hat{v} \geq 0$,

Complementary slackness:

$$\hat{u}\left(\sum_{t\in\mathcal{D}\setminus\mathcal{S}} \hat{\lambda}_t x_t + \sum_{r\in\mathcal{S}\setminus\{k\}} \hat{\lambda}_r x_r + \hat{\lambda}_k(Mx_k) - \hat{\theta}(Mx_k)\right) = 0,$$

$$\hat{v}\left(\sum_{t\in\mathcal{D}\setminus\mathcal{S}} \hat{\lambda}_t y_t + \sum_{r\in\mathcal{S}\setminus\{k\}} \hat{\lambda}_r y_r + \hat{\lambda}_k y_k - y_k\right) = 0,$$

$$\hat{w}\left(\sum_{t\in\mathcal{D}\setminus\mathcal{S}} \hat{\lambda}_t + \sum_{r\in\mathcal{S}\setminus\{k\}} \hat{\lambda}_r + \hat{\lambda}_k - 1\right) = 0,$$

$\hat{\theta}(\hat{u}(Mx_k) - 1) = 0,$

$\hat{\lambda}_r(-\hat{u}x_r + \hat{v}y_r + \hat{w}) = 0 \qquad \forall r \in \mathcal{S}\setminus\{k\},$

$\hat{\lambda}_k(-\hat{u}(Mx_k) + \hat{v}y_k + \hat{w}) = 0,$

$\hat{\lambda}_t(-\hat{u}x_t + \hat{v}y_t + \hat{w}) = 0 \qquad \forall t \in \mathcal{D}\setminus\mathcal{S}.$

All KKT constraints except $-\hat{u}x_t + \hat{v}y_t + \hat{w} \leq 0 \ \ \forall t \in \mathcal{D}\setminus\mathcal{S}$ are clearly satisfied since $\hat{\lambda}_t = 0 \ \ t \in \mathcal{D}\setminus\mathcal{S}$. Therefore, $\theta_k^{\mathcal{S}} = \theta_k^{\mathcal{D}}$ if $-\hat{u}x_t + \hat{v}y_t + \hat{w} \leq 0, \ \forall t \in \mathcal{D}\setminus\mathcal{S}$. □

From a computational geometry perspective, we interpret Result 3 as checking whether or not the points of $\mathcal{D}$ are in the same half-space $\mathcal{H}_-^{(\hat{u},\hat{v},\hat{w})}$ defined by the optimal dual solution of $P(\mathcal{S})$. Since $(x_r, y_r) \in \mathcal{H}_-^{(\hat{u},\hat{v},\hat{w})}$ $\forall r \in \mathcal{S}\backslash\{k\}$ and $(Mx_k, y_k) \in \mathcal{H}_-^{(\hat{u},\hat{v},\hat{w})}$, we need only check $(x_t, y_t)$ $\forall t \in \mathcal{D}\backslash\mathcal{S}$ as Condition (7).

We summarize the relevance of Results 1, 2, and 3 as follows.

- Results 1 and 2 show how Problems (1) and $P(\mathcal{D})$ connect, and how to determine the efficiency based on Problem (1) by using the solution of Problem $P(\mathcal{D})$. Result 3 checks whether or not a subset $\mathcal{S} \subseteq \mathcal{D}$ can represent the full set $\mathcal{D}$ (i.e. Problem $P(\mathcal{S})$ gives an identical solution to $P(\mathcal{D})$). It is the terminating condition to compute small-size LPs ($P(\mathcal{S})$'s) instead of solving a large-scale full-size LP ($P(\mathcal{D})$). Thus, we compute $P(\mathcal{S})$'s to obtain the optimal value of $P(\mathcal{D})$, which can be transformed to the optimal value of Problem (1).
- In addition to the optimal value, the optimal solutions of Problems (1) and (2), $\theta^o$, $\lambda_r^o$'s, and $(u^o, v^o, w^o)$, can be obtained from the optimal solutions of Problems $P(\mathcal{D})$ and $D(\mathcal{D})$ as shown in the proofs of Results 1 and 2. When $-u^*x_k + v^*y_k + w^* \leq 0$, we have $\theta^o = e_k^* = M\theta_k^{\mathcal{D}}$, $\lambda_r^o = \lambda_r^*$ $r \in \mathcal{D}$, and $(u^o, v^o, w^o) = (Mu^*, Mv^*, Mw^*)$. When $-u^*x_k + v^*y_k + w^* > 0$, we have $\theta^o = e_k^* = 1$, $\lambda_k^o = 1$, $\lambda_r^o = 0$ $r \in \mathcal{D}\backslash\{k\}$, and $(u^o, v^o, w^o) = \left(Mu^*, \frac{v^*}{u^*y_k + w^*}, \frac{w^*}{u^*y_k + w^*}\right)$.
- We note that Condition (7) is the *sufficient* condition for $\theta_k^{\mathcal{S}} = \theta_k^{\mathcal{D}}$, not the *necessary* condition. It is possible that Condition (7) is not satisfied (and computation continues even when reaching the optimality), even though a set $\mathcal{S}$ gives $\theta_k^{\mathcal{S}} = \theta_k^{\mathcal{D}}$. The reason is that there are alternate optimal dual solutions to $P(\mathcal{S})$, and its termination is probabilistic depending on the solution obtained. For example, an extreme efficient DMU with only itself as the reference yields degeneracy and can lead to multiple optima in the dual problem (cf. Ali, 1994; Barr and Durchholz, 1997; Cooper et al., 2007).
- In the cases we study, a linear system with $m + n + 1$ equalities provides a unique solution of $(\hat{u}, \hat{v}, \hat{w})$. According to the complementary slackness, $\hat{\lambda}_k > 0$ and $\hat{\lambda}_r > 0$ for $r \in \mathcal{S}\backslash\{k\}$ lead to $-\hat{u}(Mx_k) + \hat{v}y_k + \hat{w} = 0$ and $-\hat{u}x_r + \hat{v}y_r + \hat{w} = 0$, respectively. It follows that when their reference size is $m+n$, $(\hat{u},\hat{v},\hat{w})$ associates with a linear system with $m+n+1$ equalities and is unique. Further, when slacks exist in the optimal solution to Problem $\mathrm{P}(\mathcal{S})$, each constraint with slacks leads to zero value of the corresponding dual solution, which provides an equality contributing to solving $(\hat{u},\hat{v},\hat{w})$. As a result, a unique $(\hat{u},\hat{v},\hat{w})$ is obtained when the number of references plus the number of constraints in Problem $\mathrm{P}(\mathcal{S})$ equal $m+n$.

- Computing efficiency by solving Problem $\mathrm{P}(\mathcal{D})$ instead of (1) helps to obtain the unique optimal dual solution. Cooper et al. (2007, p.448) noted that in practice inefficient DMUs rarely have alternative optimal solutions to Problem (2) (or $\mathrm{D}(\mathcal{D})$). Typically, the projected points of the inefficient DMUs are composed with $m+n$ references or they are weakly efficient existing slacks. The unique dual optimal solution of $\mathrm{P}(\mathcal{S})$ is obtained. Therefore, when choosing a sufficiently large $M$ to formulate and solve Problem $\mathrm{P}(\mathcal{D})$, we attempt to evaluate an inefficient DMU.
- There are less than $m+n+1$ equalities in the optimal dual solution of $\mathrm{P}(\mathcal{D})$ when $(Mx_k, y_k)$'s do not project onto a $m+n-1$ dimensional facet, most likely because $(Mx_k, y_k)$ is efficient ($\theta_k^{\mathcal{D}}=1$). However, they do not require more computation time as we show in Section 5. Solving Problem $\mathrm{P}(\mathcal{D})$ instead of (1) mitigates the cases with multiple dual solutions. We note that Lovell and Rouse (2003) proposed a model similar to $\mathrm{P}(\mathcal{D})$ to fix the infeasibility of computing super efficiency (Andersen and Petersen, 1993). The super efficiency is determined in the same LP formulation (1), but the input-output data of the one evaluated is excluded from the left side of the constraints.

## 4. The algorithm: `SearchRef`

We propose a new algorithm, `SearchRef`, to solve Problem (1) that determines the efficiency of DMU $k$. Our major innovation is to consider that solving Problem (1) is a process of searching the references with the size no larger than the dimension. Our searching process is based on solving different small-size LPs instead of a full-size LP. We propose to solve a modified

LP, Problem $\mathrm{P}(\mathcal{D})$, rather than (1), to mitigate the potential of multiple dual solutions. Figure 1 is the pseudocode of `SearchRef`.

We select the data points that are the minimum of each input and the maximum of each output together with $(Mx_k, y_k)$ as the initial sample because the points are on the efficient frontier (Ali, 1993). Adding these points helps to move the initial supporting hyperplane close to the true supporting hyperplane so that it can screen out more points satisfying (and fewer points violating) Condition (7). If the number of data points selected is less than $m + n$ because some points may be superior than others in more than one input/output dimension, the empty slots are filled arbitrarily from the rest of the data set.

Sample $\mathcal{S}$ formulates Problem $\mathrm{P}(\mathcal{S})$. After solving $\mathrm{P}(\mathcal{S})$, Condition (7) gauges whether or not $\theta_k^{\mathcal{S}} = \theta_k^{\mathcal{D}}$. Failing to meet the terminating condition requires updating the members in $\mathcal{S}$, more precisely to expand $\mathcal{S}$ here. The procedure iteratively updates the samples until it selects the proper $\mathcal{S}$. Results 1 and 2 determine the optimal solutions for Problem (1) based on the $\mathrm{P}(\mathcal{S})$ solved.

We expand (update) $\mathcal{S}$ by adding $\delta \geq 1$ DMUs from $\mathcal{D} \backslash \mathcal{S}$ with the highest value of $-\hat{u}x_t + \hat{v}y_t + \hat{w}$ into the current sample. We use $-\hat{u}x_t + \hat{v}y_t + \hat{w}$ as the "infeasibility" measure for $t \in \mathcal{D} \backslash \mathcal{S}$. Violating (7) indicates that the proposed solution is infeasible with respect to $t$ (i.e. $(x_t, y_t) \in \mathcal{H}_{++}^{(\hat{u},\hat{v},\hat{w})}$). A DMU with a higher positive value implies that it is "more infeasible" with respect to Problem $\mathrm{P}(\mathcal{D})$ (or is in $\mathcal{H}_{++}^{(\hat{u},\hat{v},\hat{w})}$ and farther from $\mathcal{H}^{(\hat{u},\hat{v},\hat{w})}$); therefore, it has a higher priority of being added into the sample to fix the infeasibility.

Chen and Lai (2017) suggested selecting points that are "similar" to the improvement direction seen from the DMU evaluated as the initial sample. Although this is an effective strategy for computing Farrell's technical efficiency under the technologies with variable returns to scale, a limitation occurs when we extend it to other models and assumptions. The constant, non-decreasing, and non-increasing returns to scale assumptions allow all data points to scale either down,

Procedure `SearchRef`

// determining the efficiency of $k \in \mathcal{D}$

INPUTS: $\mathcal{D}, k, \delta, M$

OUTPUTS: $\theta_k^o, \{\lambda_r^o, r \in \mathcal{D}\}, (u^o, v^o, w^o)$

Begin

// solving P($\mathcal{D}$)

$isOpt \leftarrow$ fasle

**Selecting an initial sample $\mathcal{S} \subseteq \mathcal{D}$ such that $k \in \mathcal{S}$;**

While $isOpt =$ false Do

Solving P($\mathcal{S}$);

$\{\hat{\theta}; \hat{\lambda}_r, r \in \mathcal{S}\} \leftarrow$ optimal solution

$(\hat{u}, \hat{v}, \hat{w}) \leftarrow$ optimal dual solution

If $-\hat{u}x_t + \hat{v}y_t + \hat{w} \leq 0 \ \ \forall t \in \mathcal{D} \backslash \mathcal{S}$ Then //termination criterion

// optimal solution of P($\mathcal{D}$) obtained

$isOpt \leftarrow$ true

Else // optimal solution of P($\mathcal{D}$) NOT obtained

**Updating $\mathcal{S}$ such that $k \in \mathcal{S}$;**

End

End Do

// transforming the solution

If $-\hat{u}x_k + \hat{v}y_k + \hat{w} \leq 0$ Then

$\theta_k^o \leftarrow M\hat{\theta}$

$\lambda_r^o \leftarrow \hat{\lambda}_r, \forall r \in \mathcal{S};\ \lambda_r^o \leftarrow 0, \forall r \in \mathcal{D} \backslash \mathcal{S}$

$(u^o, v^o, w^o) \leftarrow (M\hat{u}, M\hat{v}, M\hat{w})$

Else

$\theta_k^o \leftarrow 1$

$\lambda_k^o \leftarrow 1;\ \lambda_r^o \leftarrow 0, \forall r \in \mathcal{D}\{k\}$

$(u^o, v^o, w^o) \leftarrow \left(M\hat{u}, \frac{\hat{v}}{\hat{v}y_k + \hat{w}}, \frac{\hat{w}}{\hat{v}y_k + \hat{w}}\right)$

End

End

Figure 1: The pseudocode of `SearchRef`.

up, or both. Determining the similarity requires considering all possible "virtual" points after scaling, and not just the points in the data set, which increases computation time. Therefore, we suggest a random-based initial solution, which is a standard procedure in computation algorithms, to select the initial sample. Our empirical studies described in Section 5 below demonstrate the effectiveness of the random-based rule for initial selection. Indeed, selecting an initial sample is not critical. Generally, it requires several iterations to update the sample and solve the problems to catch the true references. Using the terminating condition, Condition (7), as the screening mechanism and fixing the infeasibility by updating the sample allows us to allocate a smaller sample size for the initial selection and adjust a good improving direction.

Determining the incremental size of $\delta$ to expand the sample involves a tradeoff between the number of iterations needed and the computation loading of an individual LP. Obviously, a larger increment in expanding the sample allows us to add more "infeasible" points so that it is more likely to include references. However, a larger increment or initial sample leads to solving larger-size LPs.

$M > 1$ yields loss in numerical accuracy of the efficiency determination, i.e. a larger value of $M$ indicates a more serious accuracy loss. For example, `SearchRef` will lose 2-digit accuracy when $M = 100$ because of Result 1, and a 1-digit loss for $M = 10$, so we need a higher solution accuracy for `SearchRef` including LP solving. For example, supposing that the expected numerical accuracy of efficiencies is `1.0e-5`, we set the numerical tolerance of LP and `SearchRef` as `1.0e-6` with $M = 10$. While a sufficiently large $M$ creates an inefficient point to resolve the issue of multiple dual solutions, the tradeoffs are the costs of a large $M$ for solving individual LPs and terminating the algorithm in response to the required numerical accuracy.

Our algorithm introduces new columns in Problem $\mathrm{P}(\mathcal{S})$ (or new rows in $\mathrm{D}(\mathcal{S})$) at each iteration. We can also view it as a column (or row) generation used in the mathematical programming, although here we develop it from the computational geometry perspective.

We develop `SearchRef` to compute a particular LP, Problem (1), so that the efficiency of the corresponding DMU is determined. To compute the efficiencies of all DMUs in $\mathcal{D}$, we implement SearchRef $|\mathcal{D}|$ times for each DMU. The idea differs from the procedures proposed by Barr and Durchholz (1997) and Dulá (2011), which intend to compute the efficiencies of all DMUs. As a result, SearchRef provides the flexibility to compute the efficiencies of a subset of the full data set without any setup costs, and it also serves as a sub-procedure for other algorithms. For example, SearchRef can be used in the second stage in Barr and Durchholz (1997) and Dulá (2011) to compute the efficiencies for the cases with large-scale efficient DMUs.

## 5. Empirical implementation

### 5.1 Data and implementation

We test the simulated data sets used in Dulá (2011)[2]. Each data set contains 100,000 data points representing a peer of 100,000 DMUs. The data sets represent a broad range of data characteristics, including dimension (total number of inputs and outputs) and density (percentage of efficient DMUs). There are four dimensions: 5 (2-input 3-output), 10 (5-input 5-output), 15 (7-input 8-output), and 20 (10-input 10-output), and three densities (1%, 10%, and 25%) for a total of twelve dimension-density scenarios.

We compare SearchRef with three procedures: The naïve procedure solving the full-size LPs, the one proposed by Chen and Lai (2017), which solves LPs with a fixed size, and the state of the art BuildHull (Dulá, 2011). For Chen and Lai (2017), we arbitrarily set the number of variables for each individual LP problem as 300, which is less than the limit of most free trial software. We implement four computation procedures with a Gurobi 8.0.1 solver (Gurobi, 2018) with a Julia 0.6.3 (Bezanson et al. (2017)) interface, and perform the computations on an Intel Core i7-8700K CPU with 16 GB memory and a Windows 10 operating system.

We note that solving an LP using with the Gurobi solver involves both problem formulation

2 http://www.people.vcu.edu/~jdula/MassiveScaleDEAdata/

and problem solving. Gurobi's capability of modifying the coefficients of LPs instead of reformulating the problem helps to accelerate computation, so we implement this software feature in both `BuildHull` and the naïve procedure to obtain the computational benefit. However, when implementing `SearchRef` and Chen and Lai (2017), this feature is adopted only within computing the efficiency of a DMU, and we still need a problem formulation to solve each DMU's efficiency. To compute the efficiencies of all DMUs, `BuildHull` and the naïve procedure execute the problem formulation only once, whereas in `SearchRef` and Chen and Lai's procedure, the problem formulations are equal to the size of the data.[3] Having designed `SearchRef` to compute the efficiency for one DMU, we formulate the problem for each DMU to maintain this spirit. Finally, to truthfully reflect the real performance of each procedure, we do not introduce data preprocesses such as removing the dominated DMUs in all procedures.

The optimal value tolerance, which represents the accuracy of the efficiency determined for all procedures, is `1.0e-5`, which is sufficient for most efficiency analyses. We set $M = 10$ in `SearchRef`.[4] To maintain the same tolerance for all other procedures, the optimal tolerance of the Gurobi solver needs to be `1.0e-6`. For the same reason, the numerical tolerance for the terminating condition (7) is `1.0e-6`. We set the initial sample size as $m + n + 1$ and arbitrarily select 100 as the incremental size to expand the sample.[5]

5.2 Computation performance

Table 1 reports the computation performance for the four procedures under the twelve dimension-density scenarios. We define performance as the total computation time in seconds needed to compute the efficiencies of 100,000 DMUs. Table 1 also lists the computation time for both phases of implementing `BuildHull`. Even by implementing the feature of modifying the coefficients, the three procedures outperform the naïve procedure. Both `BuildHull` and

---

[3] Technically, it is possible to implement the problem formulation only once for `SearchRef` and for Chen and Lai's procedure and save computation time.

[4] Lovell and Rouse (2003) also provided some suggestions for selecting $M$, but extra computation is needed.

[5] The source codes are publicly available at https://github.com/wen-chih/SearchRef.jl.

Table 1: Comparing computation performance (in sec.).

| dimension | density | Full-size | Chen & Lai | BuildHull | | | SearchRef |
|---|---|---|---|---|---|---|---|
| | (%) | total | total | phase 1 | phase 2 | total | total |
| 5 | 1 | 4,890 | 1,175 | 47 | 49 | 96 | 399 |
| | 10 | 5,007 | 1,459 | 289 | 388 | 677 | 543 |
| | 25 | 4,990 | 1,526 | 759 | 823 | 1,582 | 642 |
| 10 | 1 | 7,572 | 2,041 | 61 | 70 | 131 | 387 |
| | 10 | 8,078 | 2,427 | 393 | 563 | 956 | 556 |
| | 25 | 8,231 | 2,505 | 940 | 1,200 | 2,140 | 674 |
| 15 | 1 | 10,402 | 2,828 | 73 | 84 | 157 | 489 |
| | 10 | 11,386 | 3,808 | 470 | 755 | 1,225 | 687 |
| | 25 | 11,623 | 3,712 | 1,176 | 1,516 | 2,692 | 765 |
| 20 | 1 | 15,010 | 4,182 | 101 | 114 | 215 | 674 |
| | 10 | 17,180 | 5,764 | 679 | 1,063 | 1,742 | 951 |
| | 25 | 17,656 | 5,603 | 1,687 | 2,217 | 3,904 | 1,007 |

`SearchRef` significantly outperform the procedure proposed by Chen and Lai (2017). `BuilHull` performs well in all scenarios especially for the density of 1%; it takes 96 seconds to compute the efficiencies for a data set of 100,000 DMUs with dimension 5. `BuildHull`'s computation time significantly increases as both dimension and density increase. On the other hand, `SearchRef` performs better than the other procedures in high-dimension and high-density cases. It computes >3.5 times faster than `BuildHull` when dimension=15 and density=25%. It takes 1,007 seconds, which is >3.8 times faster than `BuildHull`, to compute 100,000 efficiencies when dimension=20 and density=25%. that the results show that `SearchRef` is the most robust over data dimension and density.

Table 2 reports the average, 90th percentile, and maximum number of iterations needed for `SearchRef` under the different dimension-density scenarios. The results represent the number of iterations needed for 1.2M instances.[6] On average, the number of iterations needed to solve the efficiency of a DMU is less than 5. The 90th percentile for the twelve dimension-density scenarios is $\leq 7$, which implies that more than 90% of the instances can be solved with less than or equal to 7 iterations. The worst case among the 1.2M instances requires 10 iterations, which implies that the largest-size LPs solved have less than 1,000 variables.

We note that `SearchRef` may fail to satisfy Condition (7) and terminate the procedure, even when the true efficiency is obtained, due to the multiple dual solutions. Table 3 classifies all instances of $\mathrm{P}(\mathcal{D})$ as "eff" and "ineff" (i.e. the DMUs with $\theta_k^{\mathcal{D}} = 1$ and $\theta_k^{\mathcal{D}} < 1$, respectively) and summarizes the number of instances and the average and standard deviation and maximum number of iterations needed for "eff" and "ineff" for the twelve dimension-density scenarios. We adopt the feature in Gurobi to check the optimal solutions and slacks of the constraints. The "eff" instances are self-referenced without any slack in the optimal solution (i.e. they relate to the $\mathrm{P}(\mathcal{D})$'s with multiple optimal dual solutions). Although the "ineff" instances may have reference sizes less than $m + n$, their size of reference and number of constraints with slacks are all equal to $m + n$, which implies that a unique dual solution is obtained. Moreover, slacks exist for a very high percentage of instances, most of the slacks are associated with the outputs. Surprisingly, the "eff" instances converge faster than the "ineff" instances in terms of the average and maximum number of iterations needed. This shows that multiple dual solutions are not a computational problem for `SearchRef`.

[6] Twelve scenarios each with 100,000 LPs to solve.

Table 2: Convergence of `SearchRef`.

| | dimension = 5 | | | dimension = 10 | | | dimension = 15 | | | dimension = 20 | | |
|---|---|---|---|---|---|---|---|---|---|---|---|---|
| density | avg | 90th* | max | avg | 90th | max | avg | 90th | max | avg | 90th | max |
| 1 | 2.99 | 4 | 5 | 2.57 | 3 | 5 | 2.61 | 3 | 5 | 3.07 | 3 | 5 |
| 10 | 4.08 | 6 | 9 | 3.61 | 5 | 7 | 3.44 | 5 | 7 | 4.07 | 5 | 7 |
| 25 | 4.91 | 7 | 10 | 4.29 | 6 | 8 | 3.89 | 5 | 8 | 4.29 | 5 | 8 |

* the 90th percentile

Table 3: Mitigating the impact of multiple dual solutions.

| | | eff | | | | ineff | | | |
|---|---|---|---|---|---|---|---|---|---|
| dimension | density | no. | avg. | std. | max | avg. | std. | max | p-value* |
| 5 | 1 | 45 | 2.49 | 0.51 | 3 | 2.99 | 0.63 | 5 | 0.000 |
| | 10 | 170 | 3.02 | 0.78 | 5 | 4.08 | 1.18 | 9 | 0.000 |
| | 25 | 258 | 3.39 | 0.88 | 5 | 4.92 | 1.35 | 10 | 0.000 |
| 10 | 1 | 88 | 2.14 | 0.35 | 3 | 2.57 | 0.61 | 5 | 0.000 |
| | 10 | 266 | 2.48 | 0.54 | 4 | 3.61 | 1.07 | 7 | 0.000 |
| | 25 | 393 | 2.79 | 0.68 | 4 | 4.30 | 1.04 | 8 | 0.000 |
| 15 | 1 | 247 | 2.26 | 0.44 | 3 | 2.61 | 0.63 | 5 | 0.000 |
| | 10 | 961 | 2.70 | 0.56 | 4 | 3.45 | 1.12 | 7 | 0.000 |
| | 25 | 1,554 | 2.88 | 0.66 | 5 | 3.91 | 1.26 | 8 | 0.000 |
| 20 | 1 | 343 | 2.33 | 0.47 | 3 | 3.08 | 0.33 | 5 | 0.000 |
| | 10 | 1,366 | 2.70 | 0.54 | 4 | 4.09 | 0.71 | 7 | 0.000 |
| | 25 | 2,483 | 2.85 | 0.59 | 5 | 4.33 | 0.87 | 8 | 0.000 |

* $H_0$: no difference in the average number of iteratins

5.3 Sensitivity analyses

We are also interested in the incremental size in our sample expansion affects SearchRef's performance We test incremental sizes with 30, 50, 100, 150, and 200 points under the twelve dimension-density scenarios and record the total computation times and average number of iterations. Table 4 reports the results.

Table 4 shows that the average number of iterations needed decreases as the incremental size increases. This observation is expected, but we note that the total computation time drops first and then increases as the size increases from 30 to 200. Even though a larger expansion is more likely to include the true references and terminate the computation, we caution that reducing the number of iterations by a larger incremental size does not necessarily reduce the time; in fact it may incur a cost (i.e. needing more time to solve a larger LP). Table 4 shows that choosing the increment is the tradeoff between the chance of catching the reference (iteration needed) and the computational burden imposed.

Table 5 reports the performance when we adopt an initial sample selection which randomly selects $m + n$ DMUs. Table 5 has the same interpretation and on the same sizes 30, 50, 100, 150, and 200 as Table 4. Table 5 shows worse performance for randomized initial selection compared to SearchRef (Table 4). Even so, randomized initial selection shows superior performance (i.e.

Table 4: Computation performance with four incremental sizes.

| | | incremental size ($\delta$) | | | | | | | | | |
|---|---|---|---|---|---|---|---|---|---|---|---|
| dimension | density | 30 | | 50 | | 100 | | 150 | | 200 | |
| | (%) | time[a] | no.[b] | time | no. | time | no. | time | no. | time | no. |
| 5 | 1 | 479 | 3.59 | 439 | 3.29 | 399 | 2.99 | 385 | 2.92 | 405 | 2.86 |
| | 10 | 611 | 5.44 | 588 | 4.65 | 543 | 4.08 | 502 | 3.70 | 515 | 3.51 |
| | 25 | 733 | 6.55 | 732 | 5.74 | 642 | 4.91 | 612 | 4.43 | 623 | 4.13 |
| 10 | 1 | 489 | 3.51 | 389 | 2.87 | 387 | 2.57 | 395 | 2.37 | 367 | 2.26 |
| | 10 | 654 | 4.82 | 570 | 4.20 | 556 | 3.61 | 553 | 3.37 | 548 | 3.21 |
| | 25 | 788 | 5.86 | 707 | 5.09 | 674 | 4.29 | 636 | 3.90 | 628 | 3.62 |
| 15 | 1 | 625 | 3.59 | 529 | 2.97 | 489 | 2.61 | 502 | 2.46 | 580 | 2.38 |
| | 10 | 802 | 4.76 | 755 | 4.20 | 687 | 3.44 | 650 | 3.21 | 653 | 3.05 |
| | 25 | 935 | 5.54 | 835 | 4.73 | 765 | 3.89 | 769 | 3.62 | 863 | 3.45 |
| 20 | 1 | 800 | 3.97 | 742 | 3.47 | 674 | 3.07 | 678 | 2.94 | 637 | 2.63 |
| | 10 | 1,159 | 5.68 | 1,037 | 4.81 | 951 | 4.07 | 952 | 3.73 | 991 | 3.54 |
| | 25 | 1,206 | 6.09 | 1,118 | 5.15 | 1,007 | 4.29 | 959 | 3.96 | 979 | 3.74 |

[a] total computation time (sec.); [b] average no. of iterations needed.

around twice as fast when dimension=20 and density=25% when the incremental size is 100) compared to `BuildHull`. Moreover, like Table 4, Table 5 shows the tradeoffs in different incremental sizes.

Table 6, reports the performance when we impose the size limit of 300 variables in solving LPs as in Chen and Lai (2017). Table 6 has the same interpretation as Table 4 and Table 5. shows that `SearchRef` with size limit significantly outperforms (i.e. up to 3 times faster) the procedure in Chen and Lai (2017). `SearchRef` with size cap in the LPs shows superior performance (i.e. twice as fast when dimension=20 and density=25% when the incremental size is 50 compared to

Table 5: Computation performance (randomized intial selection).

| | | incremental size ($\delta$) | | | | | | | | | |
|---|---|---|---|---|---|---|---|---|---|---|---|
| dimension | density | 30 | | 50 | | 100 | | 150 | | 200 | |
| | (%) | time[a] | no.[b] | time | no. | time | no. | time | no. | time | no. |
| 5 | 1 | 621 | 4.76 | 496 | 4.32 | 497 | 3.94 | 501 | 3.74 | 545 | 3.65 |
| | 10 | 867 | 6.81 | 704 | 6.01 | 650 | 5.15 | 656 | 4.73 | 676 | 4.45 |
| | 25 | 959 | 7.64 | 783 | 6.78 | 739 | 5.82 | 755 | 5.34 | 783 | 5.05 |
| 10 | 1 | 691 | 4.96 | 606 | 4.29 | 525 | 3.67 | 577 | 3.39 | 538 | 3.20 |
| | 10 | 937 | 6.72 | 764 | 5.87 | 739 | 4.99 | 742 | 4.58 | 770 | 4.35 |
| | 25 | 1,019 | 7.29 | 837 | 6.37 | 811 | 5.41 | 820 | 4.96 | 848 | 4.69 |
| 15 | 1 | 831 | 5.00 | 752 | 4.34 | 763 | 3.73 | 874 | 3.46 | 816 | 3.28 |
| | 10 | 1,133 | 6.63 | 976 | 5.73 | 953 | 4.82 | 912 | 4.40 | 929 | 4.17 |
| | 25 | 1,156 | 7.17 | 1,212 | 6.24 | 1,192 | 5.29 | 1,247 | 4.84 | 1,239 | 4.56 |
| 20 | 1 | 1,037 | 5.27 | 906 | 4.49 | 849 | 3.80 | 835 | 3.46 | 829 | 3.24 |
| | 10 | 1,460 | 7.01 | 1,467 | 6.00 | 1,250 | 4.99 | 1,206 | 4.51 | 1,229 | 4.23 |
| | 25 | 1,493 | 7.39 | 1,352 | 6.36 | 1,256 | 5.34 | 1,250 | 4.88 | 1,270 | 4.58 |

[a] total computation time (sec.); [b] average no. of iterations needed.

Table 6: Computation performance (size limit: 300 variables).

| | | incremental size ($\delta$) | | | | | | | | | |
|---|---|---|---|---|---|---|---|---|---|---|---|
| dimension | density | 30 | | 50 | | 100 | | 150 | | 200 | |
| | (%) | time[a] | no.[b] | time | no. | time | no. | time | no. | time | no. |
| 5 | 1 | 418 | 3.59 | 442 | 3.29 | 361 | 2.99 | 384 | 3.12 | 392 | 3.20 |
| | 10 | 607 | 5.44 | 534 | 4.70 | 675 | 6.03 | 782 | 7.27 | 886 | 8.02 |
| | 25 | 744 | 6.55 | 715 | 6.25 | 1,056 | 9.89 | 1,176 | 11.35 | 1,283 | 11.98 |
| 10 | 1 | 530 | 3.51 | 388 | 2.87 | 369 | 2.57 | 358 | 2.38 | 352 | 2.29 |
| | 10 | 680 | 4.82 | 563 | 4.20 | 593 | 4.21 | 813 | 6.21 | 849 | 6.24 |
| | 25 | 825 | 5.86 | 686 | 5.13 | 856 | 6.52 | 1,297 | 10.24 | 1,336 | 10.16 |
| 15 | 1 | 634 | 3.59 | 577 | 2.97 | 474 | 2.61 | 470 | 2.48 | 492 | 2.42 |
| | 10 | 839 | 4.76 | 715 | 4.21 | 728 | 4.04 | 989 | 5.88 | 1,024 | 5.90 |
| | 25 | 953 | 5.54 | 914 | 4.74 | 975 | 5.64 | 1,418 | 8.72 | 1,493 | 8.92 |
| 20 | 1 | 840 | 3.97 | 703 | 3.47 | 721 | 3.07 | 663 | 2.96 | 715 | 2.66 |
| | 10 | 1,124 | 5.68 | 995 | 4.81 | 1,068 | 4.94 | 1,563 | 7.54 | 1,592 | 7.63 |
| | 25 | 1,238 | 6.09 | 1,067 | 5.18 | 1,358 | 6.56 | 2,109 | 10.30 | 2,084 | 10.24 |

[a] total computation time (sec.); [b] average no. of iterations needed.

`BuildHull`. This suggests that the innovation and performance enhancement of `SearchRef` is not based on removing the size cap on LPs.

## 6. Conclusions

This paper presents a new computation algorithm, `SearchRef`, to determine DEA efficiencies by solving small-size LPs, based on searching the corresponding references, which are a subset of the efficient DMUs with numbers no greater than the dimension. By iteratively updating the small samples and solving the corresponding small-size LPs, the proposed algorithm vastly improves large-scale DEA computation. Empirical case studies show that the proposed algorithm computes 3 times faster than the state of the art when computing large-scale, high-dimension, and high-density problems.

The proposed algorithm utilizes the development of LP algorithms, e.g., interior point method and readily available software packages, but it does not replace them. Continual advancements in LP computation should benefit the algorithm's performance. Future research will study problems such as an ill-conditioned input-output data matrix due to data range or infinite multiplier values. Although this paper has focused on solving the input-oriented radial efficiency based on variable returns to scale, future research will also compute the efficiency measures under other assumptions and requirements. Refining the proposed algorithm should also enhance computational efficiency.

The proposed algorithm provides the flexibility to serve as a sub-procedure for, or to integrate with, other methods in the literature. For example, `BuildHull` and `SearchRef` have strengths in different densities that are unobserved in advance. We suggest using `SearchRef` to compute the efficiencies of a small random sample, e.g. 1,000 out of 100,000 DMUs, to estimate the density, and then use `BuildHull` or `SearchRef` for the remaining data based on the estimated density. In our empirical implementation with dimension=20, this hybrid strategy may take an additional 6.74 seconds on average[7] when the underlying density is 1%, but save more than 2,500 seconds from misusing `BuildHull` in the 25% density cases. Integrating the proposed algorithm into other methods is well worth investigation.

**Acknowledgment**

This research is partially supported by grants from the Ministry of Science and Technology, Taiwan (NSC 101-2628-E-009-009-MY3) and (MOST 104-2410-H-009-026-MY2). The author thanks Mr. Yueh-Shan Chung for computational assistance.

[7] $(674/100{,}000) \times 1{,}000$